\documentclass[10pt]{amsart}
\usepackage{latexsym}
\usepackage{amsmath}
\usepackage{amssymb}
\usepackage{mathrsfs}
\input xy
\xyoption{all}
\input diagxy

\def\epi{\mathop{\fam0 epi}\nolimits}

\def\Orth{\mathop{\fam0 Orth}\nolimits}
\def\PSub{\mathop{\fam0 PSub}\nolimits}
\def\Sub{\mathop{\fam0 Sub}\nolimits}
\def\dom{\mathop{\fam0 dom}\nolimits}
\def\cl{\mathop{\fam0 cl}\nolimits}
\def\Pr{\mathop{\fam0 Pr}\nolimits}

\begin{document}

\title
[Nonstandard Tools for Nonsmooth Analysis]
{Nonstandard Tools for\\ Nonsmooth Analysis}

\author{S.~S. Kutateladze}
\address[]{
Sobolev Institute of Mathematics\newline
\indent 4 Koptyug Avenue\newline
\indent Novosibirsk, 630090\newline
\indent Russia}
\email{
sskut@math.nsc.ru
}
\begin{abstract}
This is an overview of  the basic tools of nonsmooth analysis which
are grounded on  nonstandard models of set theory.
By way of illustration we give a criterion for an infinitesimally
optimal path of a general discrete dynamic system.
\end{abstract}
\dedicatory{On the Occasion of the Centenary of Leonid Kantorovich}

\date{June 4, 2012}

\maketitle

\section*{Introduction}

Analysis is the technique of differentiation and integration.
Differentiation discovers trends, and integration forecasts the future
from trends. Analysis relates to the universe, reveals the glory of the
Lord, and implies equality and smoothness.

Optimization is the choice of what is most preferable.
Nonsmooth analysis is the technique of optimization
which speaks about  the humankind, reflects the diversity  of humans,
and involves inequality and obstruction.
The list of the main techniques of nonsmooth  analysis contains  subdifferential calculus
(cp.~\cite{Clarke, SubDif}).

A model within set theory is {\it nonstandard} if the membership between the
objects of the model differs from that of the originals.   In fact
the nonstandard tools of today use a couple of set-theoretic models
simultaneously. The most popular are
{\it infinitesimal analysis\/} (cp.~\cite{Bell, IBA}) and
{\it Boolean-valued analysis} (cp.~\cite{KR, Infa}).

Infinitesimal analysis provides  us with a novel understanding for
the method of indivisibles  or monadology,
synthesizing the two approaches to calculus which belong
to the inventors.

Boolean valued analysis  originated with the
famous works by Paul Cohen on the continuum hypothesis
and distinguishes itself by  the technique of  ascending and descending,
cyclic envelopes and mixings, and $B$-sets.

Calculus  reduces forecast to numbers, which is scalarization
in modern parlance. Spontaneous solutions are often labile and rarely
optimal. Thus, nonsmooth analysis deals with inequality, scalarization
and stability.  Some aspects of the latter are revealed by the tools
of nonstandard models to be discussed.

\section*{Environment for Optimization}

The best is divine---Leibniz wrote
to Samuel Clarke:\footnote{See~\cite[p.~54]{Ariew} and
cp.~\cite{Ekeland}.}

\begin{itemize}
\item[]\small
God can produce everything that is possible or whatever does not imply a
contradiction, but he wills only to produce what is the best among things
possible.
\end{itemize}

Choosing the best, we use preferences. To optimize, we use infima and
suprema for bounded sets which is practically
the {\it least upper bound property}.
So optimization needs ordered sets and primarily
boundedly complete lattices.

To operate with preferences,
we use group structure. To aggregate and scale, we use linear structure.

All these are happily provided by the {\it reals\/} $\mathbb R$, a one-dimensional
Dedekind complete vector lattice. A Dedekind  complete vector
lattice is a {\it Kantorovich space}.

Since each number is a measure of quantity,
the idea of reducing to numbers is  of a~universal importance
to mathematics. Model theory provides justification
of the {\it Kantorovich heuristic principle\/}
that the members of his spaces are  numbers as well (cp.~\cite{BooPos} and
\cite{LVK}).

Life is inconceivable without  numerous conflicting ends and interests
to be harmonized. Thus the instances appear of multiple criteria  decision
making. It is impossible as a rule
to distinguish some particular  scalar target and ignore the rest of
them. This leads to vector optimization problems, involving
order compatible with linearity.

Linear inequality implies linearity  and order.
When combined, the two produce  an ordered  vector space.
Each linear inequality in the simplest environment  of the sort
is some half-space.
Simultaneity implies many instances and so leads to the
intersections of half-spaces.  These yield polyhedra as well as
arbitrary convex sets, identifying the theory of linear inequalities with
convexity.
\cite{Harpedonaptae}

Assume that $X$~is a vector space, $E$ is an~ordered vector space,
$f:X\rightarrow E^\bullet$ is some operator,
and $C:=\dom(f)\subset X$ is a~convex set.
 A {\it vector program\/}
 $(C,f)$  is written  as follows:
  $$
  x\in C,\ \ f(x)\rightarrow \inf\!.
 $$

The standard sociological trick  includes $(C,f)$
into a~parametric family  yielding
the  {\it Legendre trasform\/} or {\it Young--Fenchel transform\/} of $f$:
$$
f^*(l):=\sup_{x\in X}{(l(x)-f(x))},
$$
with $l\in X^{\#}$  a linear functional over $X$.
The epigraph of $f^*$  is a convex subset of
$ X^{\#}$ and so $f^*$  is  convex.
Observe that $-f^*(0)$ is the value of $(C,f)$.

A convex function is locally a positively homogeneous convex
function, a {\it sublinear functional}. Recall that
$p: X\to\mathbb R$ is sublinear whenever
$$
\epi p:=
\{ (x,\ t)\in X\times\mathbb R\mid  p(x)\le t\}
$$
is a cone. Recall that a numeric function is uniquely
determined from its epigraph.

Given
$C\subset X$, put
$$
H(C):=\{(x,\ t)\in X\times\mathbb R^+ \mid x\in tC\},
$$
the  {\it H\"ormander transform\/} of $C$.
Now, $C$ is convex if and only if $H(C)$  is a~cone.
A~space with a cone is a {\it $($pre$)$ordered vector space}.

Thus, convexity and order are intrinsic to nonsmooth analysis.

\section*{Boolean Tools in Action}

Assume that $X$ is a~real vector space,
$Y$ is a~{\it Kantorovich space}. Let $\mathbb  B:=\mathbb B(Y)$ be the {\it base\/} of~$Y$, i.e., the complete Boolean algebras of
positive projections in~$Y$; and let $m(Y)$ be the  universal completion of~$Y$.
Denote by $L (X,Y)$ the space of linear operators from $X$ to~$Y$.
In case $X$ is furnished with some $Y$-seminorm on~$X$, by $L^{(m)}(X,Y)$ we mean
the {\it space of dominated operators\/} from $X$ to~$Y$.
As usual,  $\{T\le0\}:=\{x\in X \mid Tx\le0\}$; $\ker(T)=T^{-1}(0)$ for $T:X\to Y$.
Also,
$P\in \Sub(X,Y)$ means that $P$ is {\it sublinear}, while
$P\in\PSub(X,Y)$ means that $P$ is {\it polyhedral}, i.e., finitely generated.
The superscript ${}^{(m)}$ suggests domination.

{\bf Kantorovich's Theorem.}\footnote{Cp.~\cite[p.~51]{SubDif}.}
{\sl Consider the problem of finding $\mathfrak X$ satisfying

\[
\xymatrix{
  X\ar[dr]_{B} \ar[r]^{A}
                & W  \ar@{.>}[d]^{\mathfrak X}  \\
                & Y           }
\]

{(\bf 1):}
$(\exists \mathfrak X)\ {\mathfrak X}A=B \leftrightarrow {\ker(A)\subset\ker(B)}.
$

{(\bf 2):}
{\sl If $W$ is ordered by $W_+$ and $A(X)-W_+=W_+ - A(X)=W$, then}
$$
(\exists \mathfrak X\ge 0)\ {\mathfrak X}A=B \leftrightarrow \{A\le0\}\subset\{B\le 0\}.
$$
}

{\bf The Farkas Alternative.}\footnote{Cp.~\cite[Th.~1]{Trends}.}
  {\sl
Let $X$ be a~$Y$-seminormed real vector space,
with $Y$ a~Kantorovich space.
Assume  that $A_1,\dots,A_N$ and $B$ belong to~$ L^{(m)}(X,Y)$.

Then one and  only one of the following holds:

{\rm(1)} There are   $x\in X$ and $b, b' \in \mathbb B$  such that
$b'\le b$ and
$$
 b'Bx>0, bA_1 x\le0,\dots, bA_N x\le 0.
$$

{\rm(2)} There are   positive orthomorphisms $\alpha_1,\dots,\alpha_N\in\Orth(m(Y))_+$
such that
$
B=\sum\nolimits_{k=1}^N{\alpha_k A_k}.
$
}

{\bf Theorem 1.}\footnote{Cp.~\cite[Th.~1]{Polyhedral}.}
{\sl Let $X$ be a $Y$-seminormed real vector space, with $Y$ a~Kantorovich space.
Assume given some dominated operators
$A_1,\dots,A_N,  B\in L^{(m)}(X,Y)$ and elements $u_1,\dots, u_N,v\in Y$.
The following are equivalent:

{\rm(1)} For all   $b\in \mathbb B$ the  inhomogeneous operator inequality $bBx\le bv$
is a consequence of the
 consistent simultaneous inhomogeneous operator inequalities
 $bA_1 x\le bu_1,\dots, bA_N x\le bu_N$,
i.e.,
$$
\{bB\le bv\}\supset\{bA_1\le bu_1\}\cap\dots\cap\{bA_N\le bu_N\}.
$$

{\rm(2)} There are positive orthomorphisms $\alpha_1,\dots,\alpha_N\in\Orth(m(Y))$
satisfying
$$
B=\sum\limits_{k=1}^N{\alpha_k A_k};\quad v\ge \sum\limits_{k=1}^N{\alpha_k u_k}.
$$
}

\section*{Infinitesimal Tools in Action}

Leibniz wrote about his version of calculus
that ``the difference from Archime\-des
style is only in expressions which in our method are more
straightforward and more applicable to the art of
invention.''

Nonstandard analysis has the two main advantages:
it ``kills quantifiers''  and  it produces
the new notions that are impossible within a single model of set theory.
By way of example let us turn to
the nonstandard presentations of
Kuratow\-ski--Painlev\'e  limits  and the concept of infinitesimal
optimality.

Recall that the central concept of Leibniz was that of
a~{\it monad}.\footnote{Cp.~\cite{Monads}.}
In nonstandard analysis the monad $\mu(\mathscr F)$ of a standard
filter $\mathscr F$ is the intersection of all standard elements of
$\mathscr F$.

Let
$F\subset X\times Y$
be an internal correspondence from a~standard set
$X$
to a~standard set
$Y$.
Assume given a~standard filter
$\mathscr N$
on
$X$
and a~topology
$\tau$
on~$Y$.
Put
$$
\gathered
\forall \forall (F):={}^*\{y'\mid  (\forall \,x\in\mu (\mathscr N )\cap\dom (F))
(\forall \,y\approx y')(x,y)\in F\},
\\
\exists \forall (F):={}^*\{y'\mid  (\exists \,x\in\mu (\mathscr N )\cap\dom (F))
(\forall \,y\approx y')(x,y)\in F\},
\\
\forall \exists (F):={}^*\{y'\mid  (\forall \,x\in\mu (\mathscr N )\cap\dom (F))
(\exists \,y\approx y')(x,y)\in F\},
\\
\exists \exists (F):={}^*\{y'\mid  (\exists \,x\in\mu (\mathscr N )\cap\dom (F))
(\exists \,y\approx y')(x,y)\in F\},
\endgathered
$$
with
${}^*$
symbolizing standardization
and
$y\approx y'$
standing for the {\it infinite proxitity\/} between
$y$ and $y'$ in~$\tau$, i.e. $y'\in\mu(\tau(y))$.
Call
$\rm Q_1 \rm Q_2(F)$
the
$\rm Q_1\rm Q_2$-{\it limit}
of
$F$
(here
$\rm Q _k$
$(k:=1,2)$
is one of the quantifiers
$\forall$
or
$\exists$).

Assume for instance that
$F$
is a~standard correspondence  on some element of
$\mathscr N$ and look at the $\exists\exists$-limit and the
$\forall\exists$-limit.
The former is  the {\it  limit superior\/}
or {\it upper limit};
the latter is  the {\it limit  inferior\/}
or {\it lower limit\/}
of $F$ along $\mathscr N$.

{\bf Theorem~2.}\footnote{Cp.~\cite[Sect.~5.2]{Infa}.}
{\sl If   $F$ is a~standard correspondence
then
$$
\allowdisplaybreaks
\gathered
\exists \exists (F)=\bigcap\limits_{U\in\mathscr N }\
\cl
\biggl(\,
\bigcup\limits_{x\in U}\, F(x)\biggr);
\\
\forall \exists (F)=\bigcap\limits_{U\in{\ddot{\mathscr N}}}\
\cl
\biggl(\,\bigcup\limits_{x\in U}\,F(x)\biggr),
\endgathered
$$
where
$\ddot{\mathscr N}$
is the grill of a filter
$\mathscr N$ on $X$,
i.e., the family comprising all subsets of
$X$ meeting $\mu(\mathscr N)$.
}

Convexity of harpedonaptae was stable in the sense that no variation
of stakes within the surrounding rope can ever spoil the convexity
of the tract to be surveyed.

Stability is often tested by perturbation or introducing various  epsilons
in appropriate places. One of the earliest excursions in this direction
is connected with the classical  Hyers--Ulam stability theorem for
$\varepsilon$-convex functions.
Exact calculations with epsilons and sharp estimates
are often bulky and slightly mysterious.

Assume given a~convex operator
$f:X\to E^\bullet$
and a~point
$\overline x$
in the effective domain
$\dom(f):=\{x\in X\mid  f(x)<+\infty\}$
of
~$f$.
Given
$\varepsilon \ge 0$
in the positive cone
$E_+$
of
$E$,
by the
$\varepsilon $-{\it subdifferential\/}
of~$f$
at
~$\overline x$
we mean the set
$$
\partial_\varepsilon f(\overline x)
:=\big\{T\in L(X,E)\mid
(\forall x\in X)(Tx-f(x)\le T\overline x -f(\overline x)+\varepsilon) \big\}.
$$

The usual subdifferential $\partial f(\overline x)$ is the
intersection:
$$
\partial f(\overline x):=\bigcap\limits_{\varepsilon \ge 0} \partial_\varepsilon f(\overline x).
$$
In topological setting we use continuous operators, replacing
$L(X,E)$ with $\mathscr L(X,E)$.

Some cones~$K_1$
and~$K_2$
in a topological vector space~$X$
are {\it in general position\/}
provided that

{\bf (1)}~the algebraic span of $K_1$
and~$K_2$
is some subspace
$X_0\subset X$;
i.e.,
$X_0=K_1-K_2=K_2-K_1$;

{\bf (2)}~the subspace
$X_0$
is complemented; i.e., there exists a continuous projection
$P:X\rightarrow X$
such that
$P(X)=X_0$;

{\bf (3)}~$K_1$
and~$K_2$
constitute a nonoblate pair in~$X_0$.

Finally, observe that the two
nonempty convex sets $C_1$ and $C_2$  are {\it in
general position\/} if so are their H\"ormander transforms
$H(C_1)$ and $H(C_2)$.

{\bf Theorem~3.}\footnote{Cp.~\cite[Th.~4.2.8]{SubDif}.}~{\sl Let $f_1:X\times Y\rightarrow E^\bullet$ and
$f_2:Y\times Z\rightarrow E^\bullet$ be convex operators and
$\delta,\varepsilon\in E^+$. Suppose that the convolution
$f_2\vartriangle f_1$ is $\delta$-exact at some point $(x,y,z)$; i.e.,
$\delta+(f_2\vartriangle f_1)(x,y)=f_1(x,y)+f_2(y,z)$. If, moreover, the
convex sets $\epi(f_1,Z)$ and $\epi(X,f_2)$ are in general
position, then }
$$
\partial_\varepsilon(f_2\vartriangle f_1)(x,y)=
\bigcup_{\substack{\varepsilon_1\ge 0,
\varepsilon_2\ge 0,\\ \varepsilon_1+\varepsilon_2=\varepsilon+\delta}}
\partial_{\varepsilon_2}f_2(y,z)\circ\partial_{\varepsilon_1}f_1(x,y).
$$

Some alternatives
are suggested
by actual infinities, which is illustrated with the conception
of {\it infinitesimal subdifferential} and {\it infinitesimal optimality}.

Distinguish
some downward-filtered subset
~$\mathscr E$ of
$E$
that is composed of positive elements.
Assuming
$E$ and~$\mathscr E$
 standard, define the {\it monad\/}
$\mu (\mathscr E)$ of $\mathscr E$ as
$\mu (\mathscr E):=\bigcap\{[0,\varepsilon ]\mid  \varepsilon \in
{}^\circ\!\mathscr E\}$.
The members of $\mu(\mathscr E)$ are {\it positive
infinitesimals\/}  with respect to~$\mathscr E $.
As usual,
${}^\circ\!\mathscr E$
denotes the external set of~all standard members of
~$E$,
the {\it standard part\/} of
~$\mathscr E$.

Assume that the monad $\mu (\mathscr E )$
is an external cone over ${}^\circ \mathbb R $ and, moreover,
$\mu (\mathscr E)\cap{}^\circ\! E=0$.
In application, $\mathscr E $ is usually
the filter of order-units of $E$.
The relation of
{\it infinite proximity\/} or
{\it infinite closeness\/}
between the members of $E$ is introduced as follows:
$$
e_1 \approx e_2 \leftrightarrow e_1 -e_2 \in\mu
(\mathscr E )\ \&\ e_2 -e_1 \in\mu (\mathscr E ).
$$

Now
$$
Df(\overline x):=\bigcap\limits_{\varepsilon \in{}^\circ \mathscr E }\,
\partial _\varepsilon f(\overline x)=
\bigcup\limits_{\varepsilon \in\mu (\mathscr E )}\,
\partial _\varepsilon f(\overline x),
$$
which is the {\it infinitesimal subdifferential} of
$f$ at  $\overline x$.
The elements of
$Df(\overline x)$ are
{\it infinitesimal subgradients\/}
of $f$ at
~$\overline x$.

{\bf Theorem~4.}\footnote{Cp.~\cite[Th.~4.6.14]{SubDif}.}
 {\sl Let $f_1:X\times Y\rightarrow E^\bullet$ and
$f_2:Y\times Z\rightarrow E^\bullet$ be convex operators.
Suppose that the convolution
$f_2\vartriangle f_1$ is infinitesimally exact at some point $(x,y,z)$; i.e.,
$
(f_2\vartriangle f_1)(x,y)\approx f_1(x,y)+f_2(y,z).
$
If, moreover, the
convex sets $\epi(f_1,Z)$ and $\epi(X,f_2)$ are
in general
position then}
$$
D(f_2\vartriangle f_1)(x,y)=
Df_2(y,z)\circ Df_1(x,y).
$$

Assume that there exists a limited value~$e:=\inf_{x\in C}f(x)$
of some program~$(C,f)$.
A feasible point~$x_0$
is called an
{\it infinitesimal solution\/} if
$f(x_0)\approx e$,
i.e., if
$f(x_0)\leq f(x)+\varepsilon$
for every
$x\in C$
and every standard~$\varepsilon\in \mathscr E$.

{\sl  A point~$x_0\in X$
is an infinitesimal solution of the unconstrained problem~$f(x)\rightarrow\inf$
if and only if
$0\in D f(x_0)$}.

Consider some {\it Slater regular program}
$$
\Lambda x=\Lambda\bar x,\quad  g(x)\leq 0,\quad  f(x)\rightarrow\inf;
$$
i.e., first,
$\Lambda\in L(X,\mathfrak X)$
is a linear operator with values in some vector space~$\mathfrak X$,
the mappings~$f:X\rightarrow E^\bullet $
and
$g:X \rightarrow F^\bullet $
are convex operators (for the sake of convenience we assume that
$\dom (f)=\dom (g)=X$);
second,
$F$
is an Archimedean ordered vector space,
$E$
is a standard Kantorovich space of bounded elements; and, at last,
the element~$g(\bar x)$
with some feasible point~$\bar x$
is a strong order unit in~$F$.

{\bf Theorem~5.}\footnote{Cp.~\cite[Sect. 5.7]{Infa}.}
{\sl A~feasible point~$x_0$
is an infinitesimal solution of a Slater regular program
if and only if the following system of conditions is compatible:
$$
\gathered
\beta\in L^+(F,E),\quad  \gamma\in L(\mathfrak X,E),\quad  \gamma
g(x_0)\approx 0,\\
0\in D f(x_0)+ D(\beta\circ g)(x_0)+\gamma\circ\Lambda.
\endgathered
$$
}

By way of illustration look at
the general problem of optimizing
discrete dynamic systems.

Let
$X_0,\dots,X_N$
be some topological vector spaces, and
let $G_k : X_{k-1}\rightrightarrows X_k$
be a~nonempty convex correspondence
for all $k:=1,\dots,N$.
The collection~$G_1,\dots,G_N$
determines the {\it dynamic family\/} of
processes ${(G_{k,l})}_{k<l\leq N}$,
where the correspondence
$G_{k,l}: X_k\rightrightarrows X_l$
is defined as
$$G_{k,l}:=G_{k+1}\circ\dots\circ G_l\quad  \text{if}\quad  k+1<l;$$
$$G_{k,k+1}:=G_{k+1}\quad  (k:=0,1,\dots,N-1).$$
Clearly,
$G_{k,l}\circ G_{l,m}=G_{k,m}$
for all~$k<l<m\leq N$.

A~{\it path\/} or
{\it trajectory\/}
of the above family of processes is defined to be
an ordered collection of elements~$\mathfrak x:=(x_0,\dots,x_N)$
such that
$x_l\in G_{k,l}(x_k)$
for all~$k<l\leq N$.
Moreover, we say that
$x_0$
is the {\it beginning\/} of $\mathfrak x$ and
$x_N$
is the {\it ending\/} of~$\mathfrak x$.

Let
$Z$
be a topological ordered vector space.
Consider some  convex operators~$f_k:X_k\rightarrow Z$  $(k:=0,\dots,N)$
and convex sets~$S_0\subset X_0$
and
$S_N\subset X_N$.
Assume given
 a topological Kantorovich space $E$
and a monotone sublinear operator $P: Z^{N+1}\to E^\bullet$.
Given a path~$\mathfrak x:=(x_0,\dots,x_N)$,
put
$$
\mathfrak f (\mathfrak x):=(f_0(x_0), f_1(x_1)\dots, f_k(x_N)).
$$
Let $\Pr_k: Z^{N+1}\to Z$ denote the projection of $Z^{N+1}$ to the $k$th coordinate.
Then $\Pr_k (\mathfrak f(\mathfrak x))=f_k(x_k)$ for all $k:=0,\dots, N$.

Observe that $\mathfrak f$ is a convex operator from $X$ to $Z$
which is the {\it vector target\/} of the  discrete dynamic problem
under study. Assume given a monotone sublinear operator
$P: Z^{N+1}$ to~$E^\bullet$.
A~path $\mathfrak x$ is
{\it feasible\/} provided that the beginning of $\mathfrak x$ belongs to~$S_0$
and the ending of ~$\mathfrak x$, to~$S_N$.
A~path~$\mathfrak x^0:=\bigl(x^0_0,\dots,x^0_N\bigr)$
is {\it infinitesimally optimal\/} provided that
$x^0_0\in S_0$,
$x^0_N\in S_N$,
and
$P\circ\mathfrak f$
attains an infinitesimal minimum over the set of all feasible paths.
This is an instance of a~general
{\it discrete dynamic extremal problem\/} which
consists in finding a~path of a dynamic family optimal in some sense.

Introduce the sets
$$
\gathered
C_0:=S_0\times \prod^N_{k=1} X_k;\quad  C_1:=G_1\times \prod^N_{k=2} X_k;
\\
C_2:=X_0\times G_2\times\prod^N_{k=3} X_k;\dots;\quad
C_N
:=\prod^{N-2}_{k=0} X_k\times G_N;
\\
C_{N+1}:=\prod^{N-1}_{k=1} X_k\times S_N;\quad  X:=\prod^N_{k=0} X_k.
\endgathered
$$

{\bf Theorem 6.}
{\sl
Suppose that the convex sets
$$
C_0\times E^+,\dots,C_{N+1}\times E^+
$$
are in general position
as well as the sets $
X\times \epi(P)$ and $\epi(\mathfrak f)\times E$.

A feasible path~$\bigl(x^0_0,\dots,x^0_N\bigr)$
is infinitesimally optimal if and only
if the following system of conditions
is compatible:}
$$
\gathered
\alpha_k\in \mathscr L(X_k,E),\quad \beta_k\in \mathscr L^+(Z,E)
\quad (k:=0,\dots,N);
\\
\beta\in\partial(P); \ \beta_k:=\beta\circ\Pr_k;
\\
(\alpha_{k-1},\alpha_k)\in D G_k\bigl(x^0_{k-1},x^0_k\bigr)
-\{0\}
\times D(\beta_k\circ  f_k)\bigl(x^0_k\bigr)\quad  (k:=1,\dots,N);
\\
-\alpha_0\in D S_0(x_0)+D(\beta_0\circ f_0)(x_0);
\quad
 \alpha_N\in D S_N(x_N).
\endgathered
$$

{\sc Proof.} Each infinitesimally optimal path
$u:=\big(x^0_0,\dots,x^0_N\big)$
is obviously an infinitesimally optimal solution of the program
$$
v\in C_0\cap\dots\cap C_{N+1},\ \ P\circ{\mathfrak f}(v)\rightarrow \inf.
$$
By the Lagrange principle the optimal value  of this program
is the value of  some program
$$
v\in C_0\cap\dots\cap C_{N+1},\ \ {\mathfrak g(v)}\rightarrow \inf,
$$
where $\mathfrak g (v):=\beta(\mathfrak f(v))$ for all paths $v$
with $\beta\in\partial P$.  The latter has separated targets, which case
is settled  (cp.~\cite[p.~213]{Infa}).

\bibliographystyle{plain}

\end{document}